\documentclass[12pt]{article}
\usepackage{hyperref}
\usepackage{amssymb,amsmath,amsthm}
\usepackage{enumerate}

\usepackage[letterpaper,hmargin=2.70cm,vmargin=3.2cm]{geometry}

\usepackage{setspace}
\makeatletter
\renewcommand\section{\@startsection {section}{1}{\z@}%
                                   {-3.5ex \@plus -1ex \@minus -.2ex}%
                                   {2.3ex \@plus.2ex}%
                                   {\normalfont\large\bfseries}}
\makeatother

\newtheorem{theorem}{Theorem}
\newtheorem{proposition}{Proposition}
\newtheorem{lemma}{Lemma}
\newtheorem{corollary}{Corollary}
\theoremstyle{definition}
\newtheorem{definition}{Definition}
\newtheorem{remark}{Remark}
\newtheorem{example}{}
\newtheorem*{acknowledgments}{Acknowledgments}

\newcommand{\abs}[1]{\left|#1\right|}
\newcommand{\norm}[1]{\left\|#1\right\|}

\newcommand{\inner}[2]{\left\langle#1,#2\right\rangle}

\newcommand{\cH}{{\cal H}}
\newcommand{\cc}[1]{\overline{#1}}
\numberwithin{equation}{section}

\DeclareMathOperator{\im}{Im}
\DeclareMathOperator{\dom}{Dom}
\DeclareMathOperator{\Ker}{Ker}

\DeclareMathOperator{\ran}{Ran}
\DeclareMathOperator{\Sp}{Sp}
\begin{document}
\baselineskip=15 pt

\title{Applications of M.\,G. Krein's Theory of Regular Symmetric
       Operators to Sampling Theory 
       \thanks{Research partially supported by CONACYT under Project P42553­F.}}
       \author{\textbf{Luis O. Silva and Julio H. Toloza}
       \\[5mm]
       Departamento de M\'{e}todos Matem\'{a}ticos y Num\'ericos\\
       Instituto de Investigaciones en Matem\'aticas Aplicadas y en Sistemas\\
       Universidad Nacional Aut\'onoma de M\'exico\\
       C.P. 04510, M\'exico D.F.
       \\[3mm]
       \texttt{silva@leibniz.iimas.unam.mx}\\
       \texttt{jtoloza@leibniz.iimas.unam.mx}}
\date{}
\maketitle
\begin{center}
\begin{minipage}{5in}
  \centerline{{\bf Abstract}} \bigskip The classical Kramer sampling
  theorem establishes general conditions that allow the reconstruction
  of functions by mean of orthogonal sampling formulae.  One major
  task in sampling theory is to find concrete, non trivial
  realizations of this theorem.  In this paper we provide a new
  approach to this subject on the basis of the M. G. Krein's theory of
  representation of simple regular symmetric operators having
  deficiency indices $(1,1)$.  We show that the resulting sampling
  formulae have the form of Lagrange interpolation series.  We also
  characterize the space of functions reconstructible by our sampling
  formulae. Our construction allows a rigorous treatment of certain
  ideas proposed recently in quantum gravity.
\end{minipage}
\end{center}
\newpage

\section{Introduction}
\label{sec:intro}

It has been argued recently that sampling theory might be the bridge that
would allow to reconcile the continuous nature of physical fields with the
need of discretization of space-time, as required by a yet-to-be-formulated
theory of quantum gravity \cite{kempf1,kempf2,kempf3}.
This idea, which is partially developed in \cite{kempf1} for the
one-dimensional case, introduces the use of simple regular symmetric operators
with deficiency indices $(1,1)$ to obtain orthogonal sampling formulae.
It is remarkable that the class of operators under consideration in
\cite{kempf1} had been already studied in detail by M. G.  Krein
\cite{krein1,krein2,krein3} more than 60 years ago.

The conjunction of the ideas of \cite{kempf1} and Krein's theory of symmetric
operators with equal deficiency indices suggests the means of treating
sampling theory in a new and straightforward way.
By introducing this new approach, we put in a mathematically rigorous
framework the ideas in \cite{kempf1} and propose a general method for
obtaining analytic sampling formulae associated with the self-adjoint
extensions of a simple regular symmetric operator with deficiency
indices $(1,1)$.

A seminal result in sampling theory is the Whittaker-Shannon-Kotel'nikov
(WSK) sampling theorem \cite{kotelnikov,shannon,whittaker}.
This theorem states that functions that belong to a Paley-Wiener space may
be uniquely reconstructed from their values at certain discrete sets of
points.
A general approach to WSK type formulae was developed by Kramer
\cite{kramer} based on the following result.
Given a finite interval $I=[a,b]$, let $K(y,x)\in L^2(I,dy)$ for all
$x\in\mathbb{R}$.
Assume that there exists a sequence $\{x_n\}_{n\in\mathbb{Z}}$ for which
$\{K(y,x_n)\}_{n\in\mathbb{Z}}$ forms an orthogonal basis of $L^2(I,dy)$.
Let $f$ be any function of the form
\[
f(x):=\inner{K(\cdot,x)}{g(\cdot)}_{L^2(I)}
\]
for some $g\in L^2(I,dy)$, where the brackets denotes the inner product in
$L^2(I,dy)$.
Then $f$ can be reconstructed from its samples $\{f(x_n)\}_{n\in\mathbb{Z}}$
by the formula
\[
f(x)=\sum_{n=-\infty}^\infty
        \frac{\inner{K(\cdot,x)}{K(\cdot,x_n)}_{L^2(I)}}
        {\norm{K(\cdot,x_n)}^2_{L^2(I)}}f(x_n),\qquad x\in\mathbb{R}.
\]

The search for concrete realizations of the Kramer sampling theorem, which
in some cases has proven to be a difficult task, has motivated a large amount
of literature.
See, for instance, \cite{annaby1,garcia1,garcia2,garcia3,zayed1,zayed2}
and references therein.
On the basis of the approach proposed in the present work, we obtain
analytic sampling formulae for functions belonging to linear sets of
analytic functions determined by a given regular simple symmetric operator.

A central idea in Krein's theory is the fact that any simple symmetric
operator with deficiency indices $(1,1)$, acting in a certain Hilbert space
$\cH$, defines a bijective mapping from $\cH$ onto a space of scalar functions
of one complex variable having certain analytic properties.
In this space of functions one can introduce an inner product.
In the particular case when the starting point is an entire operator
(see Definition 2), the space of functions turns out to be a Hilbert space
of entire functions known as a de Branges space \cite{debranges}.

Summing up, in this paper we present an original technique in sampling
theory based on Krein's theory of regular symmetric operators.
Concretely, we obtain an analytic sampling formula valid for functions
belonging to linear sets of analytic functions associated with
symmetric operators of the kind already mentioned.  We then focus our
attention to entire operators to provide a characterization of the
corresponding spaces of entire functions.  As a byproduct, we also
provide rigorous proofs to some of the results formally obtained in
\cite{kempf1}.  In a sense this work may be considered as introductory
to the theoretical framework.  Further rigorous results related to
\cite{kempf1}, as well as other applications to sampling theory and
inverse spectral theory, will be discussed in subsequent papers.

Our approach to sampling theory allows us to reconstruct functions
that are out of the scope of certain techniques developed recently
\cite{garcia1}. It is worth remarking that, in terms of the method derived 
in this paper, the results of \cite{garcia1} correspond to considering 
only a particular self-adjoint extension of an entire operator.
We also remark that \cite{garcia1} does not fit in the
framework proposed in \cite{kempf1} for which, in contrast, our
results indeed give some rigorous justification.

This paper is organized as follows.
In Section 2, we introduce some concepts of the theory of regular symmetric
operators.
In Section 3, we state and prove a sampling formula. In Section 4, we provide
a characterization of the Hilbert spaces associated with the class of
operators under consideration, paying particular attention to the case of
entire operators.
Finally, in Section 5, we discuss some examples.

\section{Preliminaries}
\label{sec:preliminaries}

Most of the mathematical background required in this work is based on
Krein's theory of representation of symmetric operators as accounted
in the expository book \cite{R1466698}.  In this section we review
some definitions and results from the theory of symmetric operators
and introduce the notation.

Let $\cH$ denote a separable Hilbert space whose inner product
$\inner{\cdot}{\cdot}$ will be assumed anti-linear in its
first argument.

\begin{definition}
  A closed symmetric operator $A$ with domain and range in
  $\cH$ is called:
\begin{enumerate}[(a)]
\item {\em simple} if it does not have non-trivial invariant
  subspaces on which $A$ is self-adjoint,
\item {\em regular} if every point in $\mathbb{C}$ is a point
  of regular type for $A$, that is, if the operator $(A-zI)$
  has bounded inverse for every $z\in\mathbb{C}$.
\end{enumerate}
\end{definition}

By \cite[Theorem 1.2.1]{R1466698}, an operator $A$ is simple if and only if
\[
\bigcap_{z:\im z\neq 0}\ran\left(A-zI\right)=\{0\}.
\]
Moreover, as shown in \cite[Propositions 1.3.3, 1.3.5 and 1.3.6]{R1466698}, a
simple closed symmetric operator with deficiency indices $(1,1)$ is
regular if and only if some (hence every) self-adjoint extension of $A$
within $\cal H$ has only discrete spectrum of multiplicity $1$.
Also,
\begin{equation}\label{def}
\text{dim}\left[\ran\left(A-zI\right)\right]^\perp=
\text{dim}\left[\Ker\left(A^*-\overline{z}I\right)\right] = 1
\end{equation}
for all $z\in{\mathbb C}$.
For this kind of operators it is also known that, for any given
$x\in\mathbb{R}$, there is exactly one self-adjoint extension
within $\cH$ such that $x$ is in its spectrum
\cite[Proposition 1.3.5]{R1466698}.

\begin{remark}
In what follows, whenever we consider self-adjoint extensions of a symmetric
operator $A$, we will always mean the self-adjoint restrictions of $A^*$,
in other words, the self-adjoint extensions of $A$ within $\cH$ (cf. Naimark's
theory on generalized self-adjoint extensions \cite[Appendix 1]{MR1255973}).
\end{remark}

We will denote by $\text{Sym}^{(1,1)}_\text{R}({\cal H})$ the
class of regular, simple, closed symmetric operators, defined
on $\cH$, with deficiency indices $(1,1)$.

Let $A_\sharp$ be some self-adjoint extension of
$A\in\text{Sym}^{(1,1)}_\text{R}({\cal H})$.
The generalized Cayley transform is defined as
\[
\left(A_\sharp-wI\right)\left(A_\sharp-zI\right)^{-1}.
\]
for every $w\in\mathbb{C}$ and $z\in\mathbb{C}\setminus\Sp(A_\sharp)$.
This operator has several properties \cite[pag. 9]{R1466698}.
We only mention the following one:
\begin{equation}\label{cayley2}
\left(A_\sharp-wI\right)\left(A_\sharp-zI\right)^{-1}:
\Ker\left(A^*-wI\right)\to\Ker\left(A^*-zI\right)
\end{equation}
one-to-one and onto.

A complex conjugation on $\cH$ is a bijective anti-linear operator
$C:\cH\to\cH$ such that $C^2=I$ and
$\inner{C\eta}{C\varphi}=\inner{\varphi}{\eta}$ for all $\eta,\varphi\in\cH$.
A symmetric operator $A$ is said to be real with respect to a complex
conjugation $C$ if $C\dom(A)\subseteq\dom(A)$ and $CA\varphi=AC\varphi$ for
every $\varphi\in\dom(A)$.
Clearly, the condition $C\dom(A)\subseteq\dom(A)$, along with $C^2=I$,
implies that $C\dom(A)=\dom(A)$.
If moreover $A$ has deficiency indices $(1,1)$, then $A^*$ is also real with
respect to $C$ (see the proof of \cite[Corollary 2.5]{MR1627806}).

\section{Sampling formulae}
\label{sec:reg sym operators}

Let us consider an operator $A\in\text{Sym}^{(1,1)}_\text{R}({\cal H})$.
Let $A_\sharp$ be some self-adjoint extension of $A$.
Given $z_0\in{\mathbb C}\setminus\Sp(A_\sharp)$ and
$\psi_0\in\Ker\left(A^*-z_0I\right)$, define
\begin{equation}\label{psi}
\psi(z):=\left(A_\sharp-z_0I\right)\left(A_\sharp-zI\right)^{-1}\psi_0
        =\psi_0 + (z-z_0)\left(A_\sharp-zI\right)^{-1}\psi_0
\end{equation}
for every $z\in{\mathbb C}\setminus\Sp(A_\sharp)$.
This vector-valued function is analytic in the resolvent set of $A_\sharp$ and,
by (\ref{cayley2}), takes values in $\Ker\left(A^*-zI\right)$ when evaluated
at $z$.
Moreover, $\psi(z)$ has simple poles at the points of
$\Sp(A_\sharp)=\Sp_\text{disc}(A_\sharp)$.
By looking at the first equality in (\ref{psi}), it is clear that the
dependency of $\psi(z)$ on $z_0$ and $\psi_0$ is rather unessential.
Indeed, for any other $z_0'\in{\mathbb C}\setminus\Sp(A_\sharp)$, we can take
$\psi_0'=\left(A_\sharp-z_0I\right)\left(A_\sharp-z_0'I\right)^{-1}\psi_0$ in
$\Ker\left(A^*-z_0'I\right)$ and thus
$\psi(z)=\psi_0'+(z-z_0')\left(A_\sharp-zI\right)^{-1}\psi_0'$ by the first
resolvent identity.

Given $z_1\in{\mathbb C}\setminus\Sp(A_\sharp)$, let us choose
some $\mu\in\cH$ such that $\inner{\psi(\cc{z_1})}{\mu}\neq
0$.  The inner product $\inner{\psi(\cc{z})}{\mu}$ then
defines an analytic function in
$\mathbb{C}\setminus\Sp(A_\sharp)$, having zeroes at a
countable set $S_\mu$ devoid of accumulation points in
$\mathbb{C}$. The set $S_\mu$ is defined as the subset of
$\mathbb{C}$ for which $\cH$ can not be written as the direct
sum of $\ran (A-zI)$ and $\text{Span}\{\mu\}$. Following
Krein, we call the element $\mu$ a {\em gauge} \cite{krein3}.

Let us define
\begin{equation}\label{thetrick}
\xi(z):=\frac{\psi(\cc{z})}{\inner{\mu}{\psi(\cc{z})}}\,
\end{equation}
for all $z\in\mathbb{C}\setminus S_\mu$.

\begin{lemma}
  The vector-valued function $\xi(z)$ is anti-analytic in
  $\mathbb{C}\setminus S_\mu$ and has simple poles at points
  of $S_\mu$.  Moreover, it does not depend on the
  self-adjoint extension of $A$ used to define $\psi(z)$.
\end{lemma}
\begin{proof}
  The first statement holds by rather obvious reasons, notice
  only that the poles of $\inner{\psi(\cc{z})}{\mu}$ coincide
  with those of $\psi(\cc{z})$. Let us pay attention to the
  second statement. Consider two different self-adjoint extensions
  $A_\sharp$ and $A_\sharp'$. We have
\[
\psi(z)=\left(A_\sharp-z_0I\right)\left(A_\sharp-zI\right)^{-1}\psi_0\quad
	\text{ and }\quad
\psi'(z)=\left(A'_\sharp-z_0I\right)\left(A'_\sharp-zI\right)^{-1}\psi_0.
\]
Since both $\psi(z)$ and $\psi'(z)$ belong to
$\Ker\left(A^*-zI\right)$ and the dimension of this subspace
is always equal to one, it follows that $\psi'(z)=g(z)\psi(z)$
for every $z\not\in\Sp(A_\sharp)\cup\Sp(A'_\sharp)$, where
$g(z)$ is a scalar function. Inserting this identity into
(\ref{thetrick}) yields $\xi(z)=\xi'(z)$.
\end{proof}

For every $z\in{\mathbb C}\setminus S_\mu$ we have the decomposition
$
{\cal H} = \ran\left(A-zI\right)\dot{+}\,\text{Span}\left\{\mu\right\},
$
in which case every element $\varphi\in\cH$ can be written as
\[
\varphi = \left[\varphi-\widehat{\varphi}(z)\mu\right] +
	  \widehat{\varphi}(z)\mu,
\]
where $\varphi-\widehat{\varphi}(z)\mu\in\ran\left(A-zI\right)$.
A simple computation shows that
the non-orthogonal projection $\widehat{\varphi}(z)$ is given by
\begin{equation}\label{transform}
\widehat{\varphi}(z)
        :=\frac{\inner{\psi(\cc{z})}{\varphi}}{\inner{\psi(\cc{z})}{\mu}}
         =\inner{\xi(z)}{\varphi}
\end{equation}
whenever $z\in{\mathbb C}\setminus S_\mu$; it is otherwise not defined.
Indeed, because of the anti-linearity of the inner product
in its first argument, the function $\widehat{\varphi}(z)$ is analytic in
$\mathbb{C}\setminus S_\mu$ and meromorphic in $\mathbb{C}$ for every
$\varphi\in\cH$. We note that $\widehat{\mu}(z)\equiv 1$.

Let us denote the linear map
$\varphi\mapsto\widehat{\varphi}(z)$ by $\Phi_\mu$ and the linear space of
functions given by (\ref{transform}) by $\Phi_\mu{\cal H}$.
Since the operator $A$ is simple, it follows that $\Phi_\mu$ is injective and,
therefore, is an isomorphism from $\cal H$ onto $\Phi_\mu{\cal H}$
\cite[Theorem 1.2.2]{R1466698}.
Moreover, $\Phi_\mu$ transforms $A$ into the multiplication operator on
$\Phi_\mu{\cal H}$, that is,
\[
\widehat{\left(A\varphi\right)}(z)
	= \left(\Phi_\mu A\Phi_\mu^{-1}\widehat{\varphi}\right)(z)
	= z\widehat{\varphi}(z),\quad \varphi\in\dom(A).
\]

\begin{proposition}\label{interpolation}
  Assume $S_\mu\cap\mathbb{R}=\emptyset$.  Let $\{x_n\}$ be
  the spectrum of any self-adjoint extension $A_\sharp$ of
  $A$.  Then, for any analytic function $f(z)$ that belongs to
  $\Phi_\mu{\cal H}$, we have
\begin{equation}\label{sampling}
f(z) = \sum_{x_n\in\Sp(A_\sharp)}
       \frac{\inner{\xi(z)}{\xi(x_n)}}{\norm{\xi(x_n)}^2} f(x_n)\,.
\end{equation}
The convergence in (\ref{sampling}) is uniform over compact subsets
of ${\mathbb C}\setminus S_\mu$.
\end{proposition}
\begin{proof}
  Fix some arbitrary self-adjoint extension $A_\sharp$ of $A$.
  Take another self-adjoint extension $A_\sharp'\ne A_\sharp$
  to define $\psi(z)$, that is,
$\psi(z)=\left(A_\sharp'-z_0I\right)\left(A_\sharp'-zI\right)^{-1}\psi_0$.
Arrange the elements of  $\Sp\left(A_\sharp\right)$ in a
  sequence  $\{x_n\}_{n\in M}$, where $M$ is a countable
  indexing set, and let $\eta_n$ be an eigenstate of
  $A_\sharp$ corresponding to $x_n$, i.\,e.,
  $A_\sharp\eta_n=x_n\eta_n$. Since $A^*\supset A_\sharp$,
  it follows that $\eta_m\in\Ker\left(A^*-x_mI\right)$, where
  furthermore
  $\text{dim}\left[\Ker\left(A^*-x_mI\right)\right]=1$. On the
  other hand, since $x_m$ is not a pole of $\psi(z)$,
  $\psi(x_m)$ is well defined and belongs also to
  $\Ker\left(A^*-x_mI\right)$. Therefore, up to a factor,
  $\eta_m=\xi(x_m)$.

  Pick a sequence $\{M_k\}_{k\in\mathbb{N}}$ of subsets of $M$
  such that $M_k\subset M_{k+1}$ and $\cup_kM_k=M$. Consider
  any analytic function $f(z)\in\Phi_\mu{\cal H}$.  Since $\Phi_\mu$
  is injective, there exists a unique $\varphi\in{\cal H}$
  such that $\widehat{\varphi}(z)=f(z)$.  Clearly,
\[
\left|\widehat{\varphi}(z)-
\sum_{n\in M_k}\frac{\inner{\xi(z)}{\xi(x_n)}}{\norm{\xi(x_n)}^2}
                     \inner{\xi(x_n)}{\varphi}\right|
        \le \norm{\xi(z)}\norm{\varphi-\sum_{n\in M_k}\frac{1}
                     {\norm{\xi(x_n)}^2}
        \inner{\xi(x_n)}{\varphi}\xi(x_n)}.
\]
The second factor in the r.\,h.\,s. of this inequality does
not depend on $z$ and obviously tends to zero as $k\to\infty$.
Since $\xi(z)$ is continuous on ${\mathbb C}\setminus S_\mu$, the uniform
convergence of (\ref{sampling}) has been proven.
\end{proof}

\begin{remark}
Krein asserts that, for any operator in
$\text{Sym}^{(1,1)}_\text{R}({\cal H})$, one can always choose $\mu$ so that
$S_\mu\cap\mathbb{R}=\emptyset$ \cite[Theorem 8]{krein2}.
\end{remark}

Clearly, the sequence $\{x_n\}$ could be replaced by any other
sequence $\{z_n\}$ of complex numbers for which
$\left\{\xi(z_n)\right\}$ is an orthogonal basis of $\cal H$ (for a
related discussion, see \cite{garcia2}).  In our case, the question of
whether such a sequence exists or not is answered to the affirmative
by invoking the self-adjoint extensions of the operator $A$. Notice
also that, since for any real number $x$ there is a self-adjoint
extension of $A$ containing $x$ in its spectrum, it follows that every
real point can be taken as a sampling point. In contrast with this,
the construction of \cite{garcia1} has a fixed set of sampling points.

Below we show that the interpolation formula (\ref{sampling})
is indeed a Lagrange interpolation series.

\begin{proposition}
  Under the hypotheses of Proposition~\ref{interpolation},
  there exists a complex function $G(z)$, analytic in
  $\mathbb{C}\setminus S_\mu$ (hence $\mathbb{R}$) and having
  simple zeroes at $\Sp(A_\sharp)$, such that
\[
f(z) = \sum_{x_n\in\Sp(A_\sharp)}
                \frac{G(z)}{(z-x_n)G'(x_n)} f(x_n),
\]
for every $f(z)\in\Phi_\mu{\cal H}$.
\end{proposition}
\begin{proof}
Given $\{x_n\}_{n\in M}=\Sp(A_\sharp)$, set
$\psi(z)=\left(A_\sharp-z_0I\right)\left(A_\sharp-zI\right)^{-1}\psi_0$
for some $z_0$ such that $\im z_0\neq 0$ and
$\psi_0\in\Ker\left(A^*-z_0I\right)$. Define
\[
G(z):=\frac{1}{\inner{\psi(\cc{z})}{\mu}}.
\]
This function has simples zeroes at the poles of $\psi(\cc{z})$, that is,
at points of the set $\{x_n\}_{n\in M}$. Also, $\xi(z)=\cc{G(z)}\psi(\cc{z})$.
Moreover, we can write
\[
\psi(\cc{z})=\frac{\eta(\cc{z})}{\cc{z}-x_n}\quad\text{ and }\quad
G(z)=(\cc{z}-x_n)F(\cc{z}),
\]
where $\eta(z):=
(z-x_n)\left(A_\sharp-z_0I\right)\left(A_\sharp-zI\right)^{-1}\psi_0$
is analytic at $z=x_n$, $\eta(x_n)\neq 0$, and $F(x_n)=G'(x_n)$. Thus,
a straightforward computation shows that
\[
\frac{\inner{\xi(z)}{\xi(x_n)}}{\norm{\xi(x_n)}^2}
	= \frac{G(z)}{(z-x_n)G'(x_n)}
          \frac{\inner{\eta(\cc{z})}{\eta(x_n)}}{\norm{\eta(x_n)}^2}
\]
so it remains to verify that the last factor above equals one.
By the Cauchy integral formula, we have
\begin{align*}
\eta(x_n)
        &= \frac{1}{2\pi i}\oint_{|w-x_n|=\epsilon}\frac{\eta(w)}{w-x_n}\,dw\\
        &= \frac{1}{2\pi i}\oint_{|w-x_n|=\epsilon}
           \left(A_\sharp-z_0I\right)\left(A_\sharp-wI\right)^{-1}\psi_0\,dw\\
        &=-\left(A_\sharp-z_0I\right)\left(
           \frac{1}{2\pi i}\oint_{|w-x_n|=\epsilon}
           \left(wI-A_\sharp\right)^{-1}dw\right)\psi_0\\
        &=-\left(x_n-z_0\right)P_n\psi_0,
\end{align*}
where $P_n$ denotes the orthoprojector onto the eigenspace associated
with $x_n$. Therefore,
\begin{align*}
\inner{\eta(\cc{z})}{\eta(x_n)}
        &=-\left(x_n-z_0\right)\inner{\eta(\cc{z})}{P_n\psi_0}\\
        &=-\left(x_n-z_0\right)\left(z-x_n\right)
           \inner{P_n\left(A_\sharp-z_0I\right)
           \left(A_\sharp-\cc{z}I\right)^{-1}\psi_0}{\psi_0}\\
        &= \left|x_n-z_0\right|^2\inner{P_n\psi_0}{\psi_0}.
\end{align*}
Finally, it is clear that
$\inner{\eta(x_n)}{\eta(x_n)}=\left|x_n-z_0\right|^2\inner{P_n\psi_0}{\psi_0}$.
\end{proof}

\begin{remark}
Notice that the function $G(z)$ is defined up to a constant factor.
In particular, one may adjust it so that $G'(x_k)=1$, where $x_k$ is a fixed
eigenvalue of $A_\sharp$. Thus, a computation like the one in the proof above
shows that
\begin{equation}\label{gz}
G(z)=(z-x_k)\frac{\inner{\xi(z)}{\xi(x_k)}}{\norm{\xi(x_k)}^2}.
\end{equation}
This identity may be useful in some applications; see Example 2 below.
\end{remark}

\section{Spaces of analytic functions}
\label{sec:entire operators}

In this section we characterize the set of functions given by the
mapping $\Phi_\mu$. Note that in the general case $\Phi_\mu\cH$ is a
space of analytic functions with simples poles in a subset of $S_\mu$.

Let ${\cal R}\subset{\cal H}$ be the linear space of elements $\varphi$ for
which $\widehat{\varphi}(z)$ is analytic on $\mathbb{R}$. As a consequence of
\cite[Corollary 1.2.1]{R1466698}, it follows that
\begin{equation}\label{isometry}
\inner{\varphi}{\eta}
     = \int_{-\infty}^\infty\cc{\widehat{\varphi}(x)}\,\widehat{\eta}(x)\,
          dm(x)
\end{equation}
for any $\varphi,\eta\in{\cal R}$ and
$m(x)=\inner{E_x\mu}{\mu}$, where $E_x$ is any spectral
function of the operator $A$. That is, $\Phi_\mu{\cal R}$ is a
linear space of analytic functions in ${\mathbb C}\setminus
S_\mu$ such that their restriction to $\mathbb{R}$ belong to
$L^2(\mathbb{R},dm)$; in short,
\[
\left.\Phi_\mu{\cal R}\right|_\mathbb{R}\subset L^2(\mathbb{R},dm).
\]
Moreover, in this restricted sense $\Phi_\mu$ is an isometry
from $\cal R$ into $L^2(\mathbb{R},dm)$.

The following theorem is due to Krein \cite[Theorem 3]{krein2}.

\begin{theorem}[Krein]\label{krein}
  For $A\in\text{Sym}^{(1,1)}_\text{R}({\cal H})$, assume that
  $\cc{\cal R}={\cal H}$. Consider a distribution function
  $m(x)=\inner{E_x\mu}{\mu}$, where $E_x$ is a spectral
  function of $A$.  Then the map $\Phi_\mu$ generates a bijective
  isometry from $\cal H$ onto $L^2(\mathbb{R},dm)$ if and only
  if $E_x$ is orthogonal.
\end{theorem}

This theorem deserves some comments. When
$E_x$ is orthogonal it occurs that
\[
\left.\Phi_\mu\cH\right|_\mathbb{R}=L^2(\mathbb{R},dm)
\]
in the usual sense of equivalence classes. Thus, every
function in $L^2(\mathbb{R},dm)$ is, up to a set of measure
zero with respect to $m(x)$, the restriction to $\mathbb{R}$
of a unique function that is the image under $\Phi_\mu$ of one
and only one element belonging to $\cal H$. Notice that any
orthogonal spectral function of $A$ corresponds to the
spectral function of one of its self-adjoint extensions within
$\cal H$.  Since these self-adjoint extensions have only
discrete spectrum, the inner product in $L^2(\mathbb{R},dm)$,
with $m(x)=\inner{E_x\mu}{\mu}$, reduces to an expression like
\[
\int_{-\infty}^\infty\cc{f(x)}g(x)\,dm(x)=\sum_{k}c_k\cc{f(x_k)}g(x_k)
\]
whenever $E_x$ is orthogonal. That is, the equivalence classes in these
spaces are quite broad.

The following corollary is partly a straightforward consequence of
Theorem~\ref{krein}. Notice that $S_\mu\cap\mathbb{R}=\emptyset$ implies
$\cal{R}=\cH$.

\begin{corollary}\label{hachesombrero}
  Let $A\in\text{Sym}^{(1,1)}_\text{R}({\cal H})$ and choose a
  gauge $\mu$ for this operator. Assume
  that $S_\mu\cap\mathbb{R}=\emptyset$. Let $E_x$ be one of its
  orthogonal spectral functions. Then the linear space of
  functions $\widehat{\cH}_\mu:=\Phi_\mu{\cH}$, equipped with
  the inner product
\begin{equation}\label{innerprod}
\inner{f(\cdot)}{g(\cdot)}:=\int_{-\infty}^\infty\cc{f(x)}g(x)\,dm(x)
\qquad \text{where }m(x)=\inner{E_x\mu}{\mu},
\end{equation}
is a reproducing kernel Hilbert space, with reproducing kernel
$k(z,w):=\inner{\xi(z)}{\xi(w)}$.
\end{corollary}
\begin{proof}
We only verify the last statement. Given
$f(z)=\inner{\xi(z)}{\varphi}\in\widehat{\cH}_\mu$, we have
\[
\inner{k(\cdot,w)}{f(\cdot)}
	= \int_{-\infty}^\infty\cc{k(x,w)}f(x)\,dm(x)
        = \inner{\xi(w)}{\varphi}=f(w),
\]
where the second equality follows from (\ref{isometry}).
\end{proof}

Notice that, once the operator $A$ is given,
the linear space $\widehat{\cH}_\mu$ depends only on the choice of gauge $\mu$.
By Corollary~\ref{hachesombrero}, for those gauges that obeys
$S_\mu\cap\mathbb{R}=\emptyset$, $\widehat{\cH}_\mu$ may be endowed with
different Hilbert space structures, one for each orthogonal spectral function
of the operator $A$. By (\ref{isometry}), all these Hilbert spaces are
however isometrically equivalent.

Irrespective of anyone of these Hilbert space structures, $\widehat{\cH}_\mu$
possesses the following properties.

\begin{proposition}
\label{prop-cond-1}
Suppose $\mu$ such that $S_\mu\cap\mathbb{R}=\emptyset$ and consider
$\widehat{\cH}_\mu$ equipped with an inner product of the form
(\ref{innerprod}).
\begin{enumerate}[(i)]
\item Let $w$ be a non real zero of $f(z)\in\widehat\cH_\mu$.
      Then the function $g(z):=f(z)(z-\cc{w})(z-w)^{-1}$ is also in
      $\widehat\cH_\mu$ and $\norm{g(\cdot)}=\norm{f(\cdot)}$.
\item The evaluation functional, defined by $f(\cdot)\mapsto f(z)$,
      is continuous.
\end{enumerate}
\end{proposition}

\begin{proof}
(i) Suppose that $w$ is a non real zero of $f(z)$. Since
$f(z)=\inner{\xi(z)}{\varphi}$ for some $\varphi\in{\cal H}$, it follows
that $\varphi$ is orthogonal to $\psi(\cc{w})$ and therefore
$\varphi\in\ran\left(A-wI\right)$. Note that
\[
f(z)
    = \inner{\xi(z)}{\left(A-wI\right)\left(A-wI\right)^{-1}\varphi}
    = \left(z-w\right)\inner{\xi(z)}{\left(A-wI\right)^{-1}\varphi}.
\]
In these computations we have used that $\xi(z)\in\Ker\left(A^*-\cc{z}\right)$.
Moreover,
\[
\inner{\xi(z)}{\left(A-\cc{w}I\right)\left(A-wI\right)^{-1}\varphi}
    =\left(z-\cc{w}\right)\inner{\xi(z)}{\left(A-wI\right)^{-1}\varphi}
    =\frac{z-\cc{w}}{z-w}f(z).
\]
Then $f(z)\left(z-\cc{w}\right)\left(z-w\right)^{-1}\in\widehat{\cH}_\mu$.
The equality of norms follows from the fact that the Cayley transform
$\left(A_\sharp-\cc{w}I\right)\left(A_\sharp-wI\right)^{-1}$ is an isometry.

(ii) Let $f(z),g(z)\in\widehat{\cH}_\mu$. Then $f(z)=\inner{\xi(z)}{\varphi}$
and $g(z)=\inner{\xi(z)}{\eta}$ for some $\varphi,\eta\in\cH$, and furthermore
\[
\left|f(w)-g(w)\right|=\left|\inner{\xi(w)}{\varphi-\eta}\right|
        \le \norm{\xi(w)}\norm{\varphi-\eta}
        = \norm{\xi(w)}\norm{f(\cdot)-g(\cdot)}.
\]
In other words, this result follows from the fact that $\widehat{\cH}_\mu$
is a reproducing kernel Hilbert space.
\end{proof}

\begin{definition}
An operator $A\in\text{Sym}^{(1,1)}_\text{R}({\cal H})$ is called
{\em entire} if there exists a so-called {\em entire gauge} $\mu\in\cH$
such that $\widehat{\varphi}(z)$ is an entire function for every
$\varphi\in\cH$. Equivalently, $A$ is entire if
${\cal H}=\ran\left(A-zI\right)\dot{+}\,\text{Span}\left\{\mu\right\}$
for all $z\in\mathbb{C}$.
\end{definition}

Notice that if $A$ is entire, $\xi(z)$ is a vector-valued entire
function and $\widehat\cH_\mu$ is a Hilbert space of entire functions.
The spaces of reconstructible functions treated in \cite{garcia1}
corresponds to a special situation within this particular case of our
approach.
\begin{definition}\label{def:de-branges}
A Hilbert space of entire functions is called a {\em de Branges space} if,
for every $f(z)$ in that space, the following conditions holds:
\begin{enumerate}[(i)]
\item for every non real zero $w$ of $f(z)$, the function
        $f(z)(z-\cc{w})(z-w)^{-1}$ belongs to the Hilbert space
        and has the same norm as $f(z)$,
\item the function $f^*(z):=\cc{f(\cc{z})}$ belongs to the Hilbert space
        and also has the same norm as $f(z)$;
\end{enumerate}
and furthermore,
\begin{enumerate}[(i)]
\item[(iii)] for every $w:\im w\neq 0$, the linear functional
        $f(\cdot)\mapsto f(w)$  is continuous.
\end{enumerate}
\end{definition}

There is an extensive literature concerning the properties of de Branges
spaces. We refer to \cite{debranges} for more details.

It can be shown that
$\widehat{\cal H}_\mu$ is a de Branges space for certain choices of the
entire gauge $\mu$.
There are some evidence indicating that Krein noticed this fact
\cite[pag. 209]{R1466698}. Also, some hints supporting this assertion
has been given by de Branges himself
\cite{debranges1}. We however could not find any formal proof of this
statement. Thus, for the sake of completeness and for the lack of a
proper reference, we provide a proof below.

\begin{remark}
\label{rem:real-entire-gauge}
As a consequence of \cite[Lemma 2.7.1]{R1466698}, given any self-adjoint
extension $A_\sharp$ of an operator
$A\in\text{Sym}^{(1,1)}_\text{R}({\cal H})$, one can always find a complex
conjugation $C$ for which $A_\sharp$ is real.
If follows from the proof of the cited lemma that
$C\psi(z)=\psi(\cc{z})$ when $\psi(z)$ is written in terms of the real
self-adjoint extension.  Moreover, by \cite[Theorem 2.7.1]{R1466698},
the operator $A$ is also real with respect to $C$. Since by
\cite[Corollary 2.5]{MR1627806} all the self-adjoint extensions of $A$ are
real, it follows that $C\psi(z)=\psi(\cc{z})$ for every realization of
$\psi(z)$.

If furthermore $A$ is entire, then an entire gauge $\mu$ may
be chosen real, that is, $C\mu=\mu$ (see \cite[Theorem 1]{krein1} and
also \cite[Sec.~2.7.7]{R1466698}).
\end{remark}

\begin{proposition}
  Assume that an entire operator $A$ is real with respect to
  some complex conjugation $C$ and let $\mu$ be a real entire
  gauge. Then the associated Hilbert space $\widehat{\cal H}_\mu$ is a
  de Branges space.
\end{proposition}
\begin{proof}
In view of Proposition~\ref{prop-cond-1}, we only have to verify (ii).
By Remark~\ref{rem:real-entire-gauge} we know that
$C\psi(\cc{z})=\psi(z)$ thence $C\xi(z)=\xi(\cc{z})$.
Now consider any $f(z)=\inner{\xi(z)}{\varphi}$.  Clearly
$f^*(z):=\inner{\xi(z)}{C\varphi}$ also belongs to
$\widehat{\cH}_\mu$.  Furthermore,
\[
f^*(z)=\inner{\xi(z)}{C\varphi}=\cc{\inner{C\xi(z)}{\varphi}}
      =\cc{\inner{\xi(\cc{z})}{\varphi}}=\cc{f(\cc{z})}.
\]
Since $C$ is an isometry, the equality of norms follows.
\end{proof}

Notice that we only have used the simultaneous reality of the
entire operator and its entire gauge $\mu$ in showing that
$\widehat{\cH}_\mu$ obeys (ii) of
Definition~\ref{def:de-branges}.  Indeed, this condition is
also necessary.

\begin{proposition}
\label{prop-necessary}
If $\widehat{\cH}_\mu$ is a de Branges space there is a
complex conjugation $C$ with respect to which both $A$ and
$\mu$ are real.
\end{proposition}
\begin{proof}
  Let $\widehat{C}:\widehat{\cH}_\mu\to\widehat{\cH}_\mu$ be
  defined by $(\widehat{C}f)(z)=\cc{f(\cc{z})}$, for every
  $f(z)\in\widehat{\cal H}_\mu$.  Clearly, $\widehat{C}$ is a
  complex conjugation. Moreover, the multiplication operator $\widehat{A}$,
  defined on the maximal domain in $\widehat{\cH}_\mu$ by
  $(\widehat{A}f)(z)=zf(z)$, is real with respect to
  $\widehat{C}$. Now define
  $C:=\Phi_\mu^{-1}\widehat{C}\Phi_\mu$. By construction, $C$
  is a complex conjugation in $\cH_\mu$. Since
  $(\widehat{C}\widehat{\mu})(z)\equiv
  1\equiv\widehat{\mu}(z)$, it follows that $\mu$ is real with
  respect to $C$. Finally, it is not difficult to see that
  $A= \Phi_\mu^{-1}\widehat{A}\Phi_\mu$ and therefore $A$ is real with
  respect to $C$.
\end{proof}

An alternative definition of a de Branges space is given in terms of
an arbitrary entire function $e(z)$.
The de Branges space $\cH(e)$ associated to $e(z)$ is
\[
\cH(e) =\left\{f(z)\in\text{Hol}(\mathbb{C}):f(z)/e(z)\in{\cal N},
	f^*(z)/e(z)\in{\cal N},
	\int_{-\infty}^\infty\abs{f(x)/e(x)}^2dx<\infty\right\},
\]
where $\cal N$ denotes the class of analytic functions of bounded type and
non-positive mean type in the upper half plane \cite[Chap.~1]{debranges}.
From this definition, it follows that $\cH(e)$ is a reproducing kernel
Hilbert space,
whose reproducing kernel is expressed in terms of the function $e(z)$. As
shown in \cite[Chap.~2]{debranges} (see also \cite[Sec.~5]{kaltenback})
a Hilbert space of entire functions $\tilde{\cH}$ that obeys
Definition~\ref{def:de-branges}
is unitarily equivalent to a de Branges space $\cH(e)$ with
\[
e(z) = i\sqrt{\frac{\pi}{k(w_0,w_0)\im(w_0)}}(\cc{w_0}-z)k(w_0,z),
\]
where $k(w,z)$ is the reproducing kernel of $\tilde{\cH}$ and $w_0$ is any
non-real complex number such that $k(w_0,w_0)>0$.

As discussed in \cite{debranges,kaltenback}, the multiplication operator
$\widehat{A}$ with maximal domain in $\cH(e)$ is a closed, symmetric
operator with deficiency indices $(1,1)$ and domain of codimension 0 or 1.
Since in our case the multiplication operator is unitarily equivalent to
an operator belonging to $\text{Sym}^{(1,1)}_\text{R}({\cal H})$,
the codimension of $\dom(\widehat{A})$ is necessarily equal to 0,
that is, $\widehat{A}$ is densely defined. The set of self-adjoint extensions
of $\widehat{A}$ are in one-one correspondence with the set of entire functions
$s_t(z):=-\sin t\,a(z) + \cos t\,b(z)$, $t\in[0,\pi)$, where $a(z)$ and $b(z)$
are defined by the identity $e(z)= a(z) + i b(z)$.
In terms of $s_t(z)$, we have
\begin{gather*}
\dom(\widehat{A}_t)
	= \left\{ g(z)=\frac{s_t(w_0)f(z)-s_t(z)f(w_0)}{z-w_0}:
	  f(z)\in\cH(e),\text{ fixed }w_0:s_t(w_0)\neq 0\right\},\\[1mm]
\widehat{A}_tg(z) = zg(z) + f(w_0)s_t(z).
\end{gather*}
For details, see \cite[Sec.~6]{kaltenback}. Notice that the definition of
$\dom(\widehat{A}_t)$ does not depend on the choice of $w_0$, as one
can verify by resorting to the first resolvent identity.

\section{Examples}
In this section we give two elementary illustrations of the method
developed above.

\begin{example}
Consider the semi-infinite Jacobi matrix
\begin{equation}
  \label{eq:jm}
  \begin{pmatrix}
  q_1 & b_1 & 0  &  0  &  \cdots \\[1mm]
  b_1 & q_2 & b_2 & 0 & \cdots \\[1mm]
  0  &  b_2  & q_3  & b_3 &  \\
  0 & 0 & b_3 & q_4 & \ddots\\
  \vdots & \vdots &  & \ddots & \ddots
  \end{pmatrix}\,,
\end{equation}
where $b_k>0$ and $q_k\in\mathbb{R}$ for $k\in\mathbb{N}$. Fix
an orthonormal basis $\{\delta_k\}_{k\in\mathbb{N}}$ in
$\cH$. Let $J$ be the operator in $\cH$ whose matrix
representation with respect to $\{\delta_k\}_{k\in\mathbb{N}}$
is (\ref{eq:jm}). Thus, $J$ is the minimal closed operator
satisfying
\begin{equation*}
  \langle\delta_n,J\delta_n\rangle=q_n\,,\quad
 \langle\delta_{n+1},J\delta_n\rangle=\langle\delta_n,J\delta_{n+1}\rangle
  =b_n\,,\quad\forall n\in\mathbb{N}\,.
\end{equation*}
(Consult \cite[Sec.~47]{MR1255973} for a discussion on matrix
representation of unbounded symmetric operators.) It is well
known that $J$ may have either deficiency indices $(1,1)$ or
$(0,0)$ \cite[Chap.~4, Sec.~1.2]{MR0184042}. A classical result is that
if $J$ has deficiency indices $(1,1)$, then the orthogonal
polynomials of the first kind $P_k(z)$ associated with
(\ref{eq:jm}) are such that
\begin{equation*}
  \sum_{k=0}^\infty\abs{P_k(z)}^2<\infty
\end{equation*}
uniformly in any compact domain of the complex plane
\cite[Theorem 1.3.2]{MR0184042} .
Therefore, for any $z\in\mathbb{C}$,
$\pi(z)=\sum_{k=1}^\infty P_{k-1}(z)\delta_k$ is in $\cH$.
By construction, $\pi(z)$ is in the one-dimensional space $\Ker
(J^*-zI)$. It is also known that, when the deficiency indices are (1,1),
$J$ is an entire operator and $\delta_1$ is an entire gauge for
$J$ \cite[Sec.~3.1.1 and Theorem 3.1.2]{R1466698}.

Let us find $\xi(z)$ for the operator $J$. Taking into account
(\ref{thetrick}), $\inner{\delta_1}{\xi(z)}=1$ and
$\inner{\delta_1}{\pi(\cc{z})}=1$ for all $z\in\mathbb{C}$.
Then, since both $\pi(\cc{z})$ and $\xi(z)$ are in $\Ker(J^*-\cc{z}I)$ and
$\delta_1$ is entire, $\pi(\cc{z})=\xi(z)$ for all
$z\in\mathbb{C}$. Thus, for any $\varphi$ in $\cH$,
$\varphi=\sum_{k=1}^\infty\varphi_k\delta_k$, the function
$\widehat{\varphi}(z)\in\widehat\cH_{\delta_1}$ is given by
\begin{equation*}
  \widehat\varphi(z):=\langle \pi(\cc{z}),\varphi\rangle=\sum_{k=1}^\infty
    P_{k-1}(z)\varphi_k\,,\quad z\in\mathbb{C}\,.
\end{equation*}

Clearly, if $\widehat\varphi(z)\in\widehat\cH_{\delta_1}$,
$\cc{\widehat\varphi(\cc{z})}\in\widehat\cH_{\delta_1}$.
Whence, in virtue of Proposition~\ref{prop-cond-1} and,
our space $\widehat\cH_{\delta_1}$ is a de
Branges space and then, by Proposition~\ref{prop-necessary},
$\delta_1$ is real with respect to
$C=\Phi_{\delta_1}^{-1}\widehat{C}\Phi_{\delta_1}$
($\widehat{C}$ is the conjugation in $\widehat\cH_{\delta_1}$
given in the proof of Proposition~\ref{prop-necessary}).

Formula (\ref{sampling}) is written in this case as
\begin{equation*}
  \begin{split}
f(z)&=\sum_{x_n\in \Sp\left(J_\sharp\right)}
    \frac{\langle \pi(\cc{z})\pi(x_n)\rangle}{\norm{\pi(x_n)}^2} f(x_n)\\
&=\sum_{x_n\in \Sp\left(J_\sharp\right)}\frac{f(x_n)}{\norm{\pi(x_n)}^2}
\sum_{k=0}^\infty P_k(z)
  P_k(x_n)\,,\quad z\in\mathbb{C}\,,
\end{split}
\end{equation*}
where $J_\sharp$ is a certain self-adjoint extension of $J$.

In a different setting, sampling formulae obtained on the basis of
Jacobi operators have been studied before \cite{garcia00,garcia0}.
We remark that in \cite{garcia00} the interpolation formulae differ
from the ones obtained by us. In \cite{garcia0}, Jacobi operators are
treated without using M.G. Krein's theory of entire operators.
\end{example}

\begin{example}
The entire operator used here has been taken from
\cite{R1466698} and is a particular case of an example
given by Krein in \cite{krein4}.

Consider a non-decreasing bounded function $s(t)$ such that
\begin{equation*}
  s(-\infty)=0\quad\text{and}\quad s(t-0)=s(t)\,.
\end{equation*}
Fix a function defined for any $x$ in the real interval
$(-a,a)$ by
\begin{equation*}
  F(x):=\int_{-\infty}^\infty e^{ixt}ds(t)\,.
\end{equation*}
In the linear space $\widetilde{\cal L}$ of continuous
functions in $[0,a)$ vanishing in some left neighborhood of
$a$, we define a sesquilinear form as follows
\begin{equation}
\label{eq:quasiinner-product}
  (g,f):=\int_{0}^a\int_{0}^a
  F(x-t)f(x)\cc{g(t)}dxdt\,.
\end{equation}
This form is a quasi-scalar product, i.\,e., the existence of
elements $f$ in $\widetilde{\cal L}$ such that
$f\ne 0$ and nevertheless $(f,f)=0$ is not
excluded.

Denote by $\cal D$ the set of  continuously differentiable functions
$f\in\widetilde{\cal L}$ such that $f(0)=0$ and define in
$\cal D$ the differential operator $\widetilde{A}$ by the rule
$\widetilde{A}f:=if'$. It is not difficult to show that
$(g,\widetilde{A}f)=(\widetilde{A}g,f)$ and $\cal D$ is
quasi-dense in $\widetilde{\cal L}$. Now, proceeding as in
\cite[Sec.~2.8.2]{R1466698}, one defines the space $\cal{L}$
as follows
\begin{equation*}
  {\cal L}=\widetilde{\cal{L}}\setminus\widetilde{0}\,,
  \qquad\widetilde{0}=\{f\in\widetilde{\cal L}: (f,f)=0\}\,.
\end{equation*}
In $\cal{L}$ we define an inner product by
\begin{equation}
  \label{eq:inner-product}
  \inner{\eta}{\varphi}:=(g,f)\,,
\end{equation}
where $\varphi$ and $\eta$ are equivalence classes containing
$f$ and $g$, respectively. Let $\cH$ be the completion of
$\cal{L}$ and consider in it the operator $A$ such that, for
the equivalence class $\varphi$ containing $f\in\cal D$,
$A\varphi$ is the equivalence class containing
$\widetilde{A}f$. It can be shown
that $A$ is an entire operator and that
\begin{equation*}
  \widehat{\varphi}(z)=\inner{\xi(z)}{\varphi}=\int_{0}^ae^{izt}f(t)dt\,,
\end{equation*}
where $f\in\varphi$.  This identity, together with
(\ref{eq:quasiinner-product}) and (\ref{eq:inner-product}),
determines $\xi(z)$ completely \cite[Sec.~3.2.2]{R1466698}.

Notice that, in this example, the entire gauge associated with
$\xi(z)$ remains unknown. This is not however an issue since
the sampling kernel can be computed anyway by resorting to
expression (\ref{gz}).
\end{example}

\begin{acknowledgments}
We express our gratitude to M. Ballesteros for drawing our
attention to the results of \cite{kempf1}. We also thank A. Minkin for
useful remarks that led to an improved version of this work.
\end{acknowledgments}


\begin{thebibliography}{99}

\bibitem{MR0184042}
Akhiezer, N.~I.: \emph{The classical moment problem and some related questions
  in analysis}.
\newblock Hafner Publishing Co., New York, 1965.

\bibitem{MR1255973}
Akhiezer, N.~I. and Glazman, I.~M.: \emph{Theory of linear operators in
  Hilbert space}.
\newblock Dover Publications Inc., New York, 1993.

\bibitem{annaby1}
Annaby, M. H.: On sampling theory associated with the
resolvents of singular Sturm-Liouville problems.
\newblock \emph{Proc. Amer. Math. Soc.} \textbf{131} (2002), 1803--1812.

\bibitem{debranges}
de Branges, L.: {\em Hilbert spaces of entire functions}.
\newblock Prentice-Hall Inc., Englewood Cliffs, N.~J., 1968.

\bibitem{debranges1}
de Branges, L.: Some Hilbert spaces of entire functions.
\newblock {\em Proc. Amer. Math. Soc.} {\bf 10} (1959), 840--846.

\bibitem{garcia00}
Garc\'ia, A.~G. and Hern\'andez-Medina, M.~A.: The discrete
Kramer sampling theorem and indeterminate moment problems
\newblock {\em J. Comput Appl. Math.} {\bf 134} (2001), 13--22.

\bibitem{garcia0}
Garc\'ia, A.~G. and Hern\'andez-Medina, M.~A.: Discrete
Sturm-Liouville problems, Jacobi matrices and Lagrange
interpolation series.
\newblock {\em J. Math. Anal. Appl.} {\bf 280} (2003), 221--231.

\bibitem{garcia1}
Garc\'ia, A.~G. and Hern\'andez-Medina, M.~A.: Sampling theory associated
        with a symmetric operator with compact resolvent and de Branges spaces.
\newblock {\em Mediterr. J. Math.} {\bf 2} (2005), 345--356.

\bibitem{garcia2}
Garc\'ia, A.~G. and Littlejohn, L.~L.: On analytic sampling theory.
\newblock {\em J. Comput Appl. Math.} {\bf 171} (2004), 235--246.

\bibitem{garcia3}
Garc\'ia, A.~G.: Orthogonal sampling formulas: A unified approach.
\newblock {\em SIAM Review} {\bf 42} (2000), 499--512.

\bibitem{R1466698}
Gorbachuk, M.~L. and Gorbachuk, V.~I.: \emph{ M. G. Krein's lectures on entire
        operators}. Operator Theory: Advances and Applications, {\bf 97}.
\newblock Birkha\"user Verlag, Basel, 1997.

\bibitem{kaltenback}
Kaltenb\"ack, M. and Woracek, H.: Pontryagin spaces of entire functions I.
\newblock {\em Integr. Equ. Oper. Theory} {\bf 33} (1999), 34--97.

\bibitem{kempf1}
Kempf, A.: Fields with finite information density.
\newblock {\em Phys. Rev. D} {\bf 69} (2004), 124014.

\bibitem{kempf2}
Kempf, A.: Covariant information-density cutoff in curved space-time.
\newblock {\em Phys. Rev. Lett.} {\bf 92} (2004), 221301.

\bibitem{kempf3}
Kempf, A.: On the three short-distance structures which can be described by
        linear operators.
\newblock {\em Rep. Math. Phys.} {\bf 43} (1999), 171--177.

\bibitem{kotelnikov}
Kotel'nikov, V. A.: On the carrying capacity of the ``ether''
and wire in telecomunications (in Russian).
\newblock All-Union Conference on Questions of Communications.
{\em Izd. Red. Upr. Svyazi RKKA}, Moscow, 1933.

\bibitem{kramer}
Kramer, H. P.: A generalized sampling theorem.
\newblock {\em J. Math. Phys.} {\bf 38} (1959), 68--72.

\bibitem{krein1}
Krein, M. G.: On Hermitian operators with defect numbers one (in Russian).
\newblock {\em Dokl. Akad. Nauk SSSR} {\bf 43}(8) (1944), 339--342.

\bibitem{krein2}
Krein, M. G.: On Hermitian operators with defect numbers one
II (in Russian).
\newblock {\em Dokl. Akad. Nauk SSSR} {\bf 44}(4) (1944), 143--146.

\bibitem{krein3}
Krein, M. G.: On one remarkable class of Hermitian operators (in Russian).
\newblock {\em Dokl. Akad. Nauk SSSR} {\bf 44}(5) (1944), 191--195.

\bibitem{krein4}
Krein, M. G.: Fundamental aspects of the representation theory
of Hermitian operators with deficiency indices (m,m) (in Russian).
\newblock {\em Ukrain. Mat. Zh.} {\bf 2} (1949), 3--66.

\bibitem{shannon}
Shannon, C. E.: Communication in the presence of noise.
\newblock {\em Proc. IRE} {\bf 37} (1949), 10--21.

\bibitem{MR1627806}
Simon, B.: The classical moment problem as a self-adjoint finite difference
  operator.
\newblock \emph{Adv. Math.} \textbf{137}(1) (1998), 82--203.

\bibitem{whittaker}
Whittaker, J. M.: The ``Fourier'' theory of the cardinal function.
\newblock {\em Proc. Edinburgh Math. Soc.} {\bf 1} (1929),
169--176.

\bibitem{zayed1}
Zayed, A. I.: A new role of Green's function in interpolation
and sampling theory.
\newblock {\em J. Math. Anal. Appl.} {\bf 175} (1993),
222--238.

\bibitem{zayed2}
Zayed, A. I., Hinsen, G. and Butzer, P. L.: On Lagrange interpolation and
Kramer-type sampling theorems associated with Sturm-Liouville problems.
\newblock {\em SIAM J. Appl. Math.} {\bf 50} (1990),
893--909.

\end{thebibliography}
\end{document}